\documentclass[11pt,english]{amsart}

\usepackage{amsmath,amscd,amsfonts,amsthm,amssymb}

\newcommand{\shrinkmargins}[1]{
 \addtolength{\textheight}{#1\topmargin}
 \addtolength{\textheight}{#1\topmargin}
 \addtolength{\textwidth}{#1\oddsidemargin}
 \addtolength{\textwidth}{#1\evensidemargin}
 \addtolength{\topmargin}{-#1\topmargin}
 \addtolength{\oddsidemargin}{-#1\oddsidemargin}
 \addtolength{\evensidemargin}{-#1\evensidemargin}
 }

\shrinkmargins{.7}

\newcommand{\field}[1]{\mathbb{#1}}

\newcommand{\Q}{\field{Q}}

\newcommand{\Z}{\field{Z}}

\theoremstyle{plain}
\newtheorem{thm}{Theorem}[section]
\newtheorem{prop}[thm]{Proposition}
\newtheorem{cor}[thm]{Corollary}
\newtheorem{lem}[thm]{Lemma}

\theoremstyle{definition}
\newtheorem{defn}[thm]{Definition}

\newtheorem{exmp}[thm]{Example}

\theoremstyle{remark}
\newtheorem{rem}[thm]{Remark}

\title[Power map permutations on a finite group]{Power map permutations and symmetric differences
in finite groups}

\author{M\'arton Hablicsek and Guillermo Mantilla-Soler \\ }

\address{Department of Mathematics, University of Wisconsin, 480 Lincoln drive,\newline
Madison, WI 53705 USA}

\email{hablics@math.wisc.edu}

\address{Department of Mathematics, University of British Columbia, 1984 Mathematics Road,\newline 
Vancouver, BC V6T 1Z2 Canada}

\email{mantilla@math.ubc.ca}

\date{}

\begin{document}

\maketitle

\begin{abstract}

Let $G$ be a finite group. For all $a \in \Z$, such that $(a,|G|)=1$, the function $\rho_a: G \to G$ sending $g$ to $g^a$ defines 
a permutation of the elements of $G$. Motivated by a recent generalization of Zolotarev's proof of classic quadratic reciprocity, due to Duke and Hopkins, we study the signature of the permutation $\rho_a$.
By introducing the group of conjugacy equivariant maps and the symmetric difference method on groups, we exhibit an integer $d_{G}$ such that $\text{sgn}(\rho_a)=\left(\frac{d_G}{a}\right)$ for all $G$ in a large class of groups, containing all finite nilpotent and odd order groups. 
\end{abstract}
 
\section{Introduction}

Given an odd prime $p$ and an integer $m$ not divisible by $p$, the map $l \mapsto lm$ for integers $l$ defines a permutation $\rho_{m}$ on $\Z/p\Z$. The map $m \mapsto \rho_m$ can be viewed as a homomorphism from $(\Z/p\Z)^{*}$ into the symmetric group $S_p$, and composition with the signature homomorphism from $S_p$ to $\{1, -1\}$ yields a homomorphism $\sigma$ defined on $(\Z/p\Z)^{*}$.
 Note that for $k$ a generator of $(\Z/p\Z)^{*}$ we have that $\rho_k$ is a $(p-1)$-cycle, and hence $\sigma(k)=-1$. In particular we have that $\sigma$ is nontrivial. Since $(\Z/p\Z)^{*}$ is a cyclic group of even order there is a unique nontrivial homomorphism from $(\Z/p\Z)^{*}$ into $\{1,-1\}$. It follows that $\sigma(m)$ is equal to $\left(\frac{m}{p} \right)$ where $\left(\frac{\cdot}{p}\right)$ is the Legendre symbol. Using clever 
combinatorial arguments, the above characterization of the Legendre symbol, and the Chinese Reminder Theorem, Zolotarev \cite{Zolotarev}
obtains a group-theoretic proof of quadratic reciprocity. Motivated by Zolotarev's work, Duke and Hopkins define for any finite group $G$ the 
homomorphism $ \left(\frac{\cdot}{G} \right)_{DH}:(\Z/|G|\Z)^{*} \to \pm 1$ as follows: Let $C_1,...,C_m$ be the conjugacy classes of $G$. Then, for all $a \in \Z$, 
relatively prime to $|G|$, the map $\psi_a :\{C_1,...,C_m\} \to \{C_1,...,C_m\}$ sending $C_i$ to $C_{i}^{a}$ is a permutation in $S_{m}$.
Define $\left(\frac{a}{G}\right)_{DH}$ to be the signature of $(\psi_a)$. 
Duke and Hopkins define the \textit{discriminant} $D_{G}$ of $G$ as follows: $$D_G =(-1)^{s}\prod_{C} \frac{|G|}{|C|},$$ where the product runs over real conjugacy
classes $C$, and $s$ is the number of pairs of nonreal conjugacy classes. Recall that the \textit{Kronecker symbol} is the unique extension of the Legendre symbol to a symbol $\left( \frac{n}{a} \right)$ defined for any $n,a \in \Z$ characterized by the following: 

For all integers $n, a, b$ 
\begin{itemize}
\item[(i)] \[\left( \frac{n}{a} \right)\left( \frac{n}{b} \right)=\left( \frac{n}{ab} \right),\]

\item[(ii)] \[ \left( \frac{n}{2} \right) =\begin{cases} 0 & \mbox{if $n \equiv 0 \pmod{2}$} \\ 1 & \mbox{if $n \equiv 1,7 \pmod{2}$} \\ -1 & \mbox{if $n \equiv 3,5 \pmod{2},$} \end{cases} \]

\item[(iii)] \[ \left( \frac{n}{-1} \right) =\begin{cases} 1 & \mbox{if $n \geq 0 $} \\ -1 & \mbox{if $n \le 0$,} \end{cases} \]

\item[(iv)] \[ \left( \frac{n}{0} \right) =\begin{cases} 0 & \mbox{if $n \neq 1 $} \\ 1 & \mbox{if $n = 1.$} \end{cases} \]
\end{itemize}

When $a$ is not odd and positive, some authors (e.g. \cite{DH} ) define $\left( \frac{n}{a} \right)$ only when
$n \equiv 0, 1 \pmod{4}$. All the results described here are valid under any of these two conventions, but to be consistent with \cite{DH} we only use the �restricted� Kronecker symbol.
The main result
of \cite{DH} is the following:

\begin{thm}[Duke-Hopkins, 2004]\label{DH} 

For all $a \in \Z$ with $(|G|,a)=1$
$$\left(\frac{a}{G}\right)_{DH} =\left(\frac{D_G}{a}\right)$$
where $D_G$ is the discriminant of $G$.
\end{thm}
Observe that the discriminant of $G=\Z/p\Z$ is $(-1)^{\frac{p-1}{2}}p$.
By applying Theorem $\ref{DH}$ to $G=\Z/p\Z$, Duke and Hopkins obtain the statement of quadratic reciprocity.\\

Suppose now that $G$ is a finite group of order $n$ with $m$ conjugacy classes. The homomorphism $\left(\frac{\cdot}{G}\right)_{DH}$ is a generalization of 
Zolotarev's homomorphism that depends only on the conjugacy classes of $G$. The work presented in this paper explores the more direct generalization obtained by taking powers of elements. To be more explicit, if we take powers of elements instead of powers of conjugacy classes we 
obtain a permutation in $S_n$ instead of a permutation in $S_m$.
Let $\left(\frac{\cdot}{G}\right)_{el}$ be the signature of the permutation in elements. Of course, in the case of abelian groups, $\left(\frac{\cdot}{G}\right)_{el}$
and $\left(\frac{\cdot}{G}\right)_{DH}$ are the same, however for nonabelian groups the situation is quite different. Since symmetric groups, dihedral and quaternion groups of order $8$ are rational groups (i.e. every two elements of the group generating the same cyclic subgroup are conjugate), it follows that the character $\left(\frac{\cdot}{G}\right)_{DH}$ is trivial for all
of them. On the other hand a calculation shows that $\left(\frac{\cdot}{G}\right)_{el}$ is nontrivial for $S_3, S_4, D_8$ and $Q_8$.
It is natural to ask if there is a simple characterization of $\left(\frac{\cdot}{G}\right)_{el}$ analogous to the one given by Duke and Hopkins in
Theorem \ref{DH}.

For $G$ a finite group, and non-negative integer $m$, denote the number of Sylow $2$-subgroups by $n_2(G)$ and let $f_{m}(G)=\{g\in G: g^m=1\}$. Furthermore let $\epsilon(G) =\begin{cases}
1 & \text{if} \quad |G|=2n_{2}(G)\footnotemark[1] \\
0 & \text{otherwise}
\end{cases}
$.\footnotetext[1]{These are the groups $G$ of order $2n$ where $n$ is odd and all elements of even order are involutions.}

\begin{defn}
Let $G$ be a finite group. Define 

$$d_G:= (-1)^{\frac{|G|-|f_{2}(G)|}{2}}\frac{|G|^{|f_{2}(G)|}}{n_2(G)^{\epsilon(G)}}.$$

\end{defn}

\begin{rem}\label{simple}

Notice that up to square factors the formula for $d_G$ can be rewritten as follows:

\begin{itemize}
\item[$\bullet$] $(-1)^{\frac{|G|-1}{2}}|G|$ if $G$ has odd order,
\item[$\bullet$] $(-1)^{\frac{|G|-|f_{2}(G)|}{2}}\frac{|G|}{2}$ if $|G|=2n_{2}(G)$,
\item[$\bullet$] $(-1)^{\frac{|G|-|f_{2}(G)|}{2}}$ otherwise.
\end{itemize}

We opted for the definition for $d_G$, and not for the above ``simplifications", since it is a closer analog to discriminant $D_G$ defined by Duke and Hopkins. For example, for an abelian or odd order group $G$ we have that $D_G$ and $d_{G}$ are exactly the same, not only up to square factors. Moreover, as we see next, under this definition the Kronecker symbol $\left( \frac{d_G}{\cdot} \right)$ is always defined (even in the ``restricted" sense).
\end{rem}

\begin{lem}
For any finite group $G$, we have $d_G \equiv 0$ or $1 \pmod{4}$.
\end{lem}

\begin{proof}
Let $n=|G|$. If $n$ is odd, we have that $d_G = (-1)^{\frac{n-1}{2}}n\equiv 1 \pmod{4}$. If $n$ is even, $2^{f_{2}}|d_{G}$ and since $n \equiv |f_{2}(G)| \pmod{2}$ we have that $d_G \equiv 0 \pmod{4}$ in that case.
\end{proof}

The following is the main result of this paper.

\begin{thm}\label{principal}

Let $G$ be a finite group. Suppose either that $G$ is the direct product of an odd-order group and a $2$-group or a group of order $2n$ where $n$ is odd and all elements of even order are involutions.

Then, $$\left( \frac{m}{G} \right)_{el} =\left( \frac{d_G}{m} \right)$$ for all $m$ such that $(m,|G|)=1$.

\end{thm}

\begin{rem}
Note that nilpotent and odd order groups satisfy the conditions of the Theorem. 

The properties given in Theorem \ref{principal} are sufficient but not necessary. For example, 
if $G$ isomorphic to any group of order $24$ the conclusion of Theorem \ref{principal} also holds. Furthermore, among the $148$ groups of 
order no bigger than $35$ there are only four\footnote[2]{two groups of order $30$, one of order $20$ and one of order $18$ } groups not satisfying the conclusion of Theorem \ref{principal}. We wonder if it is possible to give a complete characterization of groups for which $\left( \frac{m}{G} \right)_{el}=\left( \frac{d_G}{m} \right)$ for all $m$ such that $(m,|G|)=1$. 
\end{rem}

By using the method of symmetric differences on groups (see Corollary
\ref{dgexists}) we show that for every group $G$ there exists an integer $d_G^*$ such that 
$$\left(\frac{a}{G}\right)_{el}=\left(\frac{d_G^*}{a}\right)$$
for all $a$ satisfying $(a,|G|)=1$. 

\begin{rem}
It would be interesting to see if it is possible to write a ``nice exact'' formula for $d_G^*$. We observe 
that for all groups of order less than $36$ we have that $d_G^*=n_2(G)^{\upsilon(G)}d_G$ where $\upsilon(G)$ is either $0$ or $1$. 
\end{rem}

\section{Conjugacy equivariant maps}

\begin{defn}\label{invgroup}
Let $G$ be a finite group and let $f :G \to G$ be an injective map. We say that $f$ is \emph{conjugacy equivariant} if for all
$x,g \in G$ we have that $g^{-1}f(x)g=f(g^{-1}xg)$. We denote the set of conjugacy equivariant functions of by Sym$({G})^{G}$.
\end{defn}

Let Sym$(G)$ be the group of bijections from $G$ to itself. Notice that Sym$(G)$ is naturally a $G$-group under an action for which
Sym$({G})^{G}$ is just the subgroup of Sym$(G)$ fixed by $G$. 

\begin{rem}
Note that the function $\Gamma: a \mapsto \rho_a$ is a homomorphism from $(\Z/|G|\Z)^{*}$ to a subgroup of Sym$({G})^{G}$. If $z$ is 
a nontrivial element of the center of $G$, then the element of Sym$({G})^{G}$ defined by multiplication by $z$ is not in the image of $\Gamma$. 
\end{rem}

Let $n$ be the order and $m$ be the number of conjugacy classes of $G$. Notice that there are natural homomorphisms $\psi_{\mathcal{C}}$ and $\rho_{el}$ 
from Sym$({G})^{G}$ to $S_m$ and $S_n$ respectively. With this observation in mind we define the following quadratic characters\footnote[3]{We abuse the terminology by allowing the 
trivial homomorphism to be called quadratic.} of Sym$(G)^{G}$.
$$\left(\frac{\cdot}{G}\right)_{\mathcal{C}} : \text{Sym}({G})^{G} \stackrel{\psi_{\mathcal{C}}}\rightarrow S_m \stackrel{\text{sgn}}\twoheadrightarrow S_m/A_m =\{\pm1\}$$ and 
$$\left(\frac{\cdot}{G}\right)_{El} : \text{Sym}({G})^{G} \stackrel{\rho_{el}}\rightarrow S_n \stackrel{\text{sgn}}{\twoheadrightarrow} S_n/A_n =\{\pm1\}.$$

Our next result is a key ingredient in obtaining Theorem \ref{principal}.

\begin{thm}\label{odd}
 
Let $G$ be a finite group of odd order. Then 
$$\left(\frac{\cdot}{G}\right)_{\mathcal{C}}=\left(\frac{\cdot}{G}\right)_{El}.$$
\end{thm}

\begin{proof}

Let $f \in$ Sym$({G})^{G}$ and let $\sigma =(C_1,...,C_j)$ be a cycle that appears in the cycle decomposition of 
$\psi_{\mathcal{C}}(f)$. Let $$C_{\sigma}=\bigcup_{i=1}^{j}C_i.$$ Notice that $f$ restricts to a bijection of the set $C_{\sigma}$ and if we denote this restriction by $f_{\sigma}$, then the cycle decomposition of $f_{\sigma}$ will be a subset of the cycle 
decomposition of $\rho_{el}(f)$.
It follows that $$\rho_{el}(f)=\prod_{\sigma}f_{\sigma}$$ whenever $$\psi_{\mathcal{C}}(f)=\prod_{\sigma}\sigma.$$
It suffices, therefore, to show that \[\text{sgn}(f_{\sigma})=\text{sgn}(\sigma).\]
Since the classes in the cycle $\sigma$ all have the same size $|C_1|$, we see that $|C_{\sigma}| = j|C_1|.$
Notice that all the cycles appearing in the cycle decomposition of $f_{\sigma}$ have the same length, and that length $a$ must be a multiple of $j$, say $a=jk$. Then $|C_{\sigma}|=jkr$, where $r$ is the number of disjoint cycles of $f_{\sigma}$.
Since $|C_{\sigma}|=j|C_1|$ we obtain $kr=|C_1|$, and since $|G|$ is odd we conclude that $k$ and $r$ are odd. In particular, \[\text{sgn}(f_{\sigma})=((-1)^{jk-1})^{r}=(-1)^{j-1}=\text{sgn}(\sigma).\]
\end{proof}

\begin{cor}\label{Same4Odd&Abelian}

Let $G$ be a finite group of odd order. Then the signature of the power map permutation on conjugacy classes 
agrees with the signature of the power map permutation on elements, or in other words,
$$\left(\frac{\cdot}{G}\right)_{DH}=\left(\frac{\cdot}{G}\right)_{el}.$$
\end{cor}

\begin{proof}
Let $a \in (\Z/|G|\Z)^{*}$, and let $\rho_{a}$ be the bijection on $G$ given by raising to the $a^{\text{th}}$ power. Then $$\left(\frac{a}{G}\right)_{DH}=
\left(\frac{\rho_{a}}{G}\right)_{\mathcal{C}}=\left(\frac{\rho_{a}}{G}\right)_{El}=\left(\frac{a}{G}\right)_{el}.$$ 

\end{proof}

\section{Abelian groups, odd order groups and 2-groups}

In this section we prove Theorem \ref{principal} for abelian groups, groups of odd order, $2$-groups, and groups such that the set of elements of odd order form a subgroup of index $2$.

Let $n >1$ be an integer, and let $a$ be an integer relatively prime to $n$. We denote the multiplicative order of $a$, as an element of 
$(\Z/n\Z)^{*}$; by $o_n(a)$. Also for a given group $G$ we denote the total number of elements in $G$ of order $n$ by $G(n)$.

\begin{prop}\label{abodd} Let $G$ be a finite group which is either abelian or has odd order. 
Then, $\left( \frac{a}{G} \right)_{el} =\left( \frac{d_G}{a} \right)$ for all $(a,|G|)=1$.

\end{prop}

\begin{proof}

We have \[\left( \frac{\cdot}{G} \right)_{el}=\left( \frac{\cdot}{G} \right)_{DH}=\left( \frac{D_G}{\cdot} \right),\]
where the first equality holds by Corollary \ref{Same4Odd&Abelian} if $|G|$ is odd and it is obvious if $G$ is abelian, and the second equality holds by Theorem \ref{DH}.
It suffices, therefore, to show that $D_G = d_G$ if $|G|$ is odd or $G$ is abelian. By the above equalities we have that the sign of $D_G$ is equal to $\left( \frac{-1}{G} \right)_{el}$. Now, raising 
to the power $-1$ is a permutation of $G$ that can be written as a product of $(|G|-|f_{2}(G)|)/2$ disjoint transpositions, i.e. $\prod\limits_{ g\neq f_{2}(G)}(g,g^{-1})$.
Therefore if $G$ has odd order we have that $(-1)^{s}= (-1)^{(|G|-1)/2}$. On the other hand, odd order groups are characterized by having a unique real
class, so $D_G=(-1)^{(|G|-1)/2}|G|$. Then by definition we have that $d_G=(-1)^{(|G|-1)/2}|G|$ for $G$ an odd order group, hence 
$D_G=d_G$ for such a $G$.
If $G$ is abelian, then a real conjugacy class consists only of one element in $f_{2}(G)$. 
So $D_G=d_G$ in this case as well. 
 
\end{proof}

\begin{cor}\label{semidi}
 
Let $G$ be a finite group of order $2n$ where $n$ is odd and all elements of even order are involutions. Then, 
$\left( \frac{a}{G} \right)_{el} =\left( \frac{d_G}{a} \right)$ for all $(a,|G|)=1$.
\end{cor}

\begin{proof}

Notice first that $n_2(G)=n$, otherwise any element of order $2$ would be properly contained in its centralizer yielding to elements of even order 
that are not involutions. Let $H$ be a subgroup of $G$ of order $n$. Since $n_2(G)=n$, we have that $|f_{2}(G)|=n+1$. In particular $G$ is the 
disjoint union of $H$ and $f_{2}(G) \setminus \{1\}$. Since raising to an odd power fixes $f_{2}(G)$ elementwise, we have that
$\left( \frac{a}{G} \right)_{el}=\left( \frac{a}{H} \right)_{el}$ for all $(a,|G|)=1$. Then, by Proposition \ref{abodd} we have that
$\left( \frac{a}{G} \right)_{el}=\left( \frac{d_H}{a} \right)$ for all $(a,|G|)=1$. Our result follows since 
$$d_G= d_H(n^{(n-1)/2}2^{(n+1)/2})^2.$$

\end{proof}

\begin{prop}\label{explicit}
Let $G$ be a finite group and let $a$ be an integer coprime to the order of $G$. Then
 
$$\left(\frac{a}{G}\right)_{el}=\prod_{ \stackrel{d | |G|}{d \neq 1,2}}\left( (-1)^{\frac{\phi(d)}{o_{d}(a)}} \right)^{\frac{G(d)}{\phi(d)}}$$ 
where $\phi$ denotes the Euler totient function.

\end{prop}

\begin{proof}
Let $d$ be a divisor of $|G|$, and let $G_d$ be the set of elements of $G$ of order $d$. Notice that $G_d$ is invariant 
under the permutation $\rho_{a}$ of $G$ defined by raising to the $a^{\text{th}}$ power. Since $G$ is the disjoint union of $G_d$, where $d$ runs over 
all divisors of $|G|$,
we have that $$\left(\frac{a}{G}\right)_{el} =\prod_{d | |G|}\text{sign}(\rho_{a}\upharpoonright_{G_d}).$$ 
Moreover, we may assume that $d\neq 1,2$ since for those cases $\rho_{a}\upharpoonright_{G_d}$ is the trivial permutation. Let $x$ be an element of order $d$, and notice that $$(x,x^a,...,x^{a^{(o_{d}(a))}})$$ is the element of the cycle decomposition of
$\rho_{a}\upharpoonright_{G_d}$ that contains $x$. Since $x$ is an arbitrary element we have that $\rho_{a}\upharpoonright_{G_d}$ is a product of $\frac{G(d)}{o_{d}(a)}$ cycles
with the same cycle structure as the one containing $x$. Therefore, 
$$\text{sign}(\rho_{a}\upharpoonright_{G_d})=\left((-1)^{\left(o_{n}(a)-1\right)\left(\frac{G(d)}{o_{d}(a)}\right)} \right).$$
Since $\phi(d) | G(d)$ for all positive integer $d$, and $\phi(d)$ is even for all $d\ge 3$ the result follows. 

\end{proof}

Now we will deal with the case when $G$ is a 2-group:

\begin{prop}\label{even}
Let $G$ be a finite 2-group. Then $\left( \frac{a}{G} \right)_{el} =\left( \frac{d_G}{a} \right)$ for all $(a,|G|)=1$.

\end{prop}

\begin{proof}

Suppose that $|G|=2^n$. By Proposition \ref{abodd} we may assume that $n \geq 3$. Also, as we saw during the proof of Proposition \ref{abodd}, the permutation $\rho_{-1}$ can be written as a product of $(|G|-|f_{2}(G)|)/2$ disjoint transpositions. It follows from the definition of $d_G$ for 2-groups that \[\left( \frac{-1}{G} \right)_{el} =\left(-1\right)^{\frac{|G|-|f_{2}(G)|}{2}}=\left( \frac{d_G}{-1} \right).\]
In particular, since the Kronecker symbol is multiplicative, it is enough to prove the result whenever $a$ is an odd prime $p$. Since $f_{2}(G)$ is even, we must show that 
$$\left(\frac{p}{G}\right)_{el}=\left( \frac{-1}{p} \right)^{\frac{|G|-|f_{2}(G)|}{2}}.$$
First we show that \[\left(\frac{p}{G}\right)_{el}=\left( \frac{-1}{p} \right)^{\frac{G(4)}{2}}.\]
Let $d=2^{k}$ be a divisor of $|G|$, and suppose that $3 \leq k$. Since 
\[\left(\field{Z}/2^k\field{Z}\right)^* \cong \Z/2\Z\times \Z/2^{k-2}\Z\] we have that $\frac{\phi(d)}{o_{d}(p)}$ is even. Therefore, by 
Proposition \ref{explicit}, we have that 
\[\left(\frac{p}{G}\right)_{el}=\left(\left(-1\right)^{\frac{2}{o_4(p)}}\right)^{\frac{G(4)}{2}}=\left( \frac{-1}{p} \right)^{\frac{G(4)}{2}}.\]
Since 4 divides $|G|$, Frobenius' Theorem (see \cite{IR}) yields that 4 divides $|f_4(G)|$. Therefore $|G| \equiv |f_{4}(G)| \pmod{4}$. Furthermore $|f_{4}(G)|-|f_{2}(G)|= G(4)$ hence $|G|-|f_2(G)|\equiv G(4) \pmod{4}$ and the statement follows. \end{proof}

\section{Symmetric differences of groups}

In this section we develop the method of symmetric differences of groups, and apply it to obtain the final step in the proof of Theorem \ref{principal}. 

Let $\Delta$ be the usual symmetric difference operator in sets (i.e., $X\Delta Y=(X\cup Y)\setminus (X\cap Y)$), and let us denote the order of an element $g$ by $o(g)$. 

Next we prove that for every finite group $G$ there exist an integer $t$, and cyclic subgroups $C_1$, $C_2$, ..., $C_t$ of $G$ such that $G=\Delta_{i=1}^t C_i$. We call such decomposition $\Delta_{i=1}^t C_i$ \textit{a cyclic decomposition}. Furthermore, we say that a cyclic decomposition $\Delta_{i=1}^t C_i$ is \textit{reduced} if $C_{i} \neq C_{j}$ whenever $i\neq j$.

Notice that a finite group $G$ admits different cyclic decompositions, however: 

\begin{lem} Let $G$ be a finite group.
Then $G$ admits a unique reduced cyclic decomposition up to the ordering of the cyclic subgroups. 
\end{lem}

\begin{proof}

Let us begin with an observation. If $x$ and $y$ are two elements of a group $G$ generating the same cyclic subgroup then for every subgroup $H$ we have $x\in H$ if and only if $y\in H$. Therefore if $G_1$, ..., $G_l$ are some subgroups of $G$ then $x\in \Delta_{i=1}^l G_i$ if and only if $y\in \Delta_{i=1}^l G_i$. The main consequence of this observation is the following: let $x$ be a generator for a cyclic subgroup $C$ and let $G_1$, ..., $G_l$ be some subgroups, then $x\in \Delta_{i=1}^l G_i$ if and only if any other generator of $C$ will be contained in $\Delta_{i=1}^l G_i$. 

Now let us prove first the existence. Enumerate the cyclic subgroups of $G$ decreasingly with respect to their order. In other words, write the cyclic subgroups of $G$ as $C_1$, $C_2$, ..., $C_k=1$
where $|C_i|\geq |C_j|$ whenever $i\leq j$. Also let us choose one generator for every cyclic subgroup, so we have $k$ elements of $G$: $g_1$, ..., $g_k$. From the enumeration we can see that $o(g_i)\geq o(g_j)$ whenever $i\leq j$. 

Now construct subsets $U_{i} \subseteq G$ recursively, as follows: Set $U_{1} = C_1$ and for $i > 1$, let

\[U_{i}=\begin{cases} U_{i-1} \Delta C_i & \mbox{if } g_i\not \in U_{i-1}\\
U_{i-1} & \mbox{otherwise}\\
\end{cases}\]

We must show that $U_k =G$, so we prove that every element $g$ of $G$ is contained in $U_k$. Let $g \in G$, and let $C_i$ be the group generated by $g$. By the above observation we only have to check that $g_i\in U_k$. Clearly, $g_i \in U_i$, so if we suppose that $g_i\not\in U_k$, there is some smallest subscript $j > i$ such that $g_i\not \in U_j$. Then $g_i \in U_{j-1}$, so $U_{j-1}\not = U_j$, and thus $U_j = U_{j-1}\Delta C_j$. Therefore $g_i\in C_j$, but since $g_i$ is a generator for $C_i$, hence $C_i\subset C_j$. Thus $|C_i|<|C_j|$, which contradicts the ordering of the cyclic subgroups.

Now we prove the uniqueness. Suppose $\mathcal{B}$ and $\mathcal{C}$ are collections of distinct cyclic subgroups of $G$ such that $\Delta \mathcal{B}=G=\Delta \mathcal{C}$. If $\mathcal{B}\not=\mathcal{C}$ then choose a cyclic subgroup $B$ of largest possible order such that $B$ is contained in exactly one of the two collections. We can assume that $B \in \mathcal{B}$. Let $B=\left\langle b\right\rangle$. Now let us consider cyclic subgroups of $G$ other than $B$ containing $b$. Every such subgroup is strictly larger than $B$, and hence it lies either in both of the collections or neither. Since $b\in B$, it follows that $b$ is contained in exactly one more member of $\mathcal{B}$ than $\mathcal{C}$. This is a contradiction, since $b$ lies in an odd number of members of $\mathcal{B}$ and an odd number of members of $\mathcal{C}$.
\end{proof}

It is slightly more convenient to work with elements instead of the subgroups generated by them, therefore instead of writing $G=\Delta_{i\in I} C_i$ from now on we will write $G=\Delta_{i=1}^t \left\langle g_i\right\rangle$. With this notation:

\begin{exmp} $Q=\left\langle i\right\rangle\Delta \left\langle j\right\rangle\Delta \left\langle k\right\rangle$, where $Q$ denotes the quaternion group with the usual notations.
\end{exmp} 
\begin{exmp} $S_3=\left\langle (123)\right\rangle \Delta \left\langle (12)(3)\right\rangle \Delta \left\langle (13)(2)\right\rangle \Delta \left\langle (1)(23)\right\rangle \Delta \left\langle (1)(2)(3)\right\rangle$.
\end{exmp}

We can make the following observations:

\begin{lem}{\label{ni}}
Let $G$ be an arbitrary finite group and let us express $G$ as a symmetric difference of some cyclic subgroups, i.e. $G=\Delta_{i=1}^t \left\langle g_i\right\rangle$. Then:
\begin{itemize}

\item[(i)] $t$ is odd,

\item[(ii)] If $G$ is a nontrivial $2$-group and the decomposition is the reduced cyclic decomposition then none of the $g_i$ is equal to the identity.
\end{itemize}

\end{lem}

\begin{proof} (i) follows immediately since $1 \in \left\langle g_i\right\rangle$ for all $i$. To show (ii), notice that every non-trivial cyclic subgroup of a nontrivial finite $2$-group $G$ has exactly one involution. Furthermore, if $G=\Delta_{i=1}^t \left\langle g_i\right\rangle$, then each involution lies in an odd number of members of $\left\langle g_i\right\rangle$. Since the number of involutions in $G$ is odd, hence in the decomposiiton we need odd number of nontrivial cyclic subgroup. Furthermore, by part (i) we know that $t$ is odd thus the statement follows.\end{proof}

Furthermore, symmetric differences work nicely with direct products:

\begin{lem}
Assume that $G=H\times K$ where the orders of $H$ and $K$ are coprime. Let $H=\Delta_{i=1}^h \left\langle h_i\right\rangle$ and $K=\Delta_{j=1}^k\left\langle k_j\right\rangle$. Then $G=\Delta_{i=1...h,j=1...k} \left\langle (h_i,k_j)\right\rangle$.
\end{lem}

\begin{proof}
Since the order of $H$ and the order of $K$ are coprime we see that $(h,k)$ is in $\left\langle(h_i,k_j)\right\rangle$ if and only if $h\in \left\langle h_i\right\rangle$ and $k\in \left\langle k_j\right\rangle$. We know that both $h\in \left\langle h_i\right\rangle$ and $k\in \left\langle k_j\right\rangle$ hold for odd number of indices therefore $(h,k)\in \left\langle(h_i,k_j)\right\rangle$ holds an odd number of times. It follows that $G=\Delta_{i=1...h,j=1...k} \left\langle (h_i,k_j)\right\rangle$.
\end{proof}

We now apply symmetric differences for the remaining cases of Theorem \ref{principal}. We will restrict the power map permutation to the cyclic 
subgroups in the decomposition above.

\begin{defn}
$X$ is a power-invariant set in $G$ if for every $x\in X$ and for every $a$ which is relatively prime to the order of $G$, we have $x^a\in X$.
\end{defn}

For example, any subgroup is a power-invariant subset. Furthermore:

\begin{lem}\label{pow}
If $X$ and $Y$ are power-invariant sets in a finite group $G$, then, $X^c=G\setminus X$, $X\cap Y$, $X\cup Y$, $X\setminus Y$ and $X\Delta Y$ are power-invariant sets.
\end{lem}

\begin{proof}
The first three statements are immediate. The fourth follows from $X\setminus Y=X\cap Y^c$. Finally, the fifth comes from the third and fourth.
\end{proof}

Let $G$ be a finite group with a power-invariant subset $X$. If $a \in \Z$ 
is relatively prime to $|G|$, denote the sign of the permutation induced by $a$ on $X$ by $\left(\frac{a}{X}\right)_{el}$.

\begin{lem}\label{disj} If $X$ and $Y$ are two power-invariant sets, then:
$$\left(\frac{a}{X\Delta Y}\right)_{el}=\left(\frac{a}{X}\right)_{el}\left(\frac{a}{Y}\right)_{el}.$$
\end{lem}

\begin{proof} First we can see that if $X$ and $Y$ are disjoint power-invariant sets then $X\Delta Y=X\cup Y$. In this case a decomposition of $X\cup Y$ as a union of permutation cycles gives us a decomposition of $X$ and of $Y$ as a union of permutation cycles therefore $\left(\frac{a}{X\Delta Y}\right)_{el}=\left(\frac{a}{X\cup Y}\right)_{el}=\left(\frac{a}{X}\right)_{el}\left(\frac{a}{Y}\right)_{el}.$

In the general case by definition we know that $X\Delta Y=(X\setminus Y)\cup (Y\setminus X)$ furthermore the latter two sets are disjoint thus:
$$\left(\frac{a}{X\Delta Y}\right)_{el}=\left(\frac{a}{(X\setminus Y)\cup (Y\setminus X)}\right)_{el}=\left(\frac{a}{X\setminus Y}\right)_{el}\left(\frac{a}{Y\setminus X}\right)_{el}.$$
By the same argument the above is equal to
$$\left(\frac{a}{X}\right)_{el}\left(\frac{a}{Y}\right)_{el}\left(\frac{a}{X\cap Y}\right)_{el}^{-2}=\left(\frac{a}{X}\right)_{el} \left(\frac{a}{Y}\right)_{el}.$$

\end{proof}

As an immediate consequence we can use symmetric differences to calculate $\left(\frac{a}{G}\right)_{el}$ in terms of some cyclic subgroups of $G$:

\begin{cor}\label{symdiffandsign}
If $G=\Delta_{i=1}^t \left\langle g_i\right\rangle$, then $\left(\frac{a}{G}\right)_{el}=\prod_{i=1}^t \left(\frac{a}{\left\langle g_i\right\rangle}\right)_{el}.$
\end{cor}

Furthermore we obtain:

\begin{cor}\label{dgexists}
For every group $G$ there exists an integer $d_G^*$ such that $$\left(\frac{a}{G}\right)_{el}=\left(\frac{d_G^*}{a}\right)$$
for all $a$ satisfying $(a,|G|)=1$.
\end{cor}

\begin{proof}
By the previous lemmas we can write $G$ as $\Delta_{i=1}^t \left\langle g_i\right\rangle$ and therefore $\left(\frac{a}{G}\right)_{el}=\prod_{i=1}^t \left(\frac{a}{\left\langle g_i\right\rangle}\right)_{el}$. By Proposition \ref{abodd} it follows that 
$$\left(\frac{a}{G}\right)_{el}=\prod_{i=1}^t \left(\frac{d_{\left\langle g_i\right\rangle}}{a}\right)_{el}=\left(\frac{\prod_{i=1}^t d_{\left\langle g_i\right\rangle}}{a}\right)$$
so the corollary is proved.
\end{proof}

\begin{prop}\label{final} Assume that $G=H \times K$ where $H$ is a nontrivial finite $2$-group and $K$ is a finite group of odd order. Then, $\left( \frac{a}{G} \right)_{el} =\left( \frac{d_G}{a} \right)$ for all $(a,|G|)=1$.

\end{prop}

\begin{proof}
Let us write $H=\Delta_{i=1}^h \left\langle h_i\right\rangle$ and $K=\Delta_{j=1}^k\left\langle k_j\right\rangle$. By Lemma 
\ref{ni} we can write $H=\Delta_{i=1}^h \left\langle h_i\right\rangle$ with none of the $h_i$ equal to identity.
Then, by our previous lemmas, we have $\left(\frac{a}{G}\right)_{el}=\prod_{i=1..h,j=1...k} \left(\frac{a}{\left\langle (h_i,k_j)\right\rangle}\right)_{el}$. 
Since $h_i$ and $k_j$ have relatively prime orders, the cyclic group $\left\langle (h_i,k_j)\right\rangle$ has order $o(h_i)o(k_j)$. 
By Proposition \ref{abodd} we then have that
$$\left(\frac{a}{\left\langle (h_i,k_j)\right\rangle}\right)_{el}=\left(\frac{(-1)^{\frac{o(h_i)o(k_j)-2}{2}}}{a}\right).$$
It follows from Corollary \ref{symdiffandsign} that
$$\left(\frac{a}{G}\right)_{el}=\prod_{i=1}^h\prod_{j=1}^k\left(\frac{(-1)^{\frac{o(h_i)o(k_j)-2}{2}}}{a}\right)=\prod_{i=1}^h\left(\frac{\prod_{j=1}^k(-1)^{\frac{o(h_i)o(k_j)-2}{2}}}{a}\right).$$
Since $k$ is odd, all of the $o(k_j)$ are odd and all of the $o(h_i)$ are even, so we have that 
$$\prod_{j=1}^k(-1)^{\frac{o(h_i)o(k_j)-2}{2}}=(-1)^{\frac{o(h_i)-2}{2}}.$$ Therefore,
$$\left(\frac{a}{G}\right)_{el}=\prod_{i=1}^h\left(\frac{(-1)^{\frac{o(h_i)-2}{2}}}{a}\right)=\left(\frac{a}{H}\right)_{el}.$$
On the other hand, $H$ is a nontrivial 2-group therefore, by Proposition \ref{even}
$$\left(\frac{a}{H}\right)_{el}=\left(\frac{d_H}{a}\right)=\left(\frac{(-1)^{\frac{|H|-|f_2(H)|}{2}}}{a}\right)=\left(\frac{(-1)^{\frac{|K||H|-|f_2(H)|}{2}}}{a}\right),$$
since $K$ has odd order. Furthermore, since $G=H\times K$ hence $|f_2(G)|=|f_2(H)|>1$ and $|K||H|=|G|$, thus
$$\left(\frac{a}{H}\right)_{el}=\left(\frac{(-1)^{\frac{|G|-|f_{2}(G)|}{2}}|G|^{|f_{2}(G)|}}{a}\right)=\left(\frac{d_G}{a}\right).$$ In other words $\left( \frac{a}{G} \right)_{el} =\left( \frac{d_G}{a} \right)$.
\end{proof}
Theorem \ref{principal} now follows by Proposition \ref{final} and
Corollary \ref{semidi}.

\section{Further remarks}

Let $n$ be a positive integer. Notice that by identifying $(\Z/n\Z)^{*}$ with the Galois group $\text{Gal}(
\Q(\zeta_{n})/\Q)$ any non-trivial homomorphism from $(\Z/n\Z)^{*}$ to $\{\pm1\}$ corresponds an isomorphism $$\text{Gal}(
(\Q(\sqrt{d})/\Q) \to \{\pm1\},$$ where $\Q(\sqrt{d}) \subseteq \Q(\zeta_{n})$ is
a quadratic subextension. Therefore, by Chebotarev's density theorem (see \cite{marcus}), any homomorphism from $(\Z/n\Z)^{*}$ to $\{\pm1\}$ is of the form
$\left( \frac{d}{\cdot} \right)$ for some integer $d$. In particular to any finite group $G$, by considering the homomorphisms, \[ \left(\frac{\cdot}{G} \right)_{DH}:(\Z/|G|\Z)^{*} \to \pm 1\quad \mbox{and} \quad \left(\frac{\cdot}{G} \right)_{el}:(\Z/|G|\Z)^{*} \to \pm 1,\] we can associate a pair number fields $K_{G}$ and $L_{G}$ of degree at most two. In this context the main result of Duke and Hopkins is the following: Let $\Delta_{G}$ be the determinant of the character table of $G$. Then, \[K_{G} =\Q(\Delta_{G})=\Q(\sqrt{D_{G}}).\] It follows that up to square factors $D_G$ is the discriminant of $K_{G}$, and hence the name discriminant of $G$. Notice that by Corollary \ref{dgexists} we have that \[L_{G} =\Q(\sqrt{d^{*}_{G}}).\] However, it seems that in practice the integer $d^{*}_{G}$ is not so easy to calculate. On the other hand for a big class of finite groups, namely the ones described in the hypothesis of Theorem \ref{principal}, we have that \[L_{G} =\Q(\sqrt{d_{G}})\] and as we observed in Remark \ref{simple} the integer $d_{G}$ is quite simple to calculate. We wonder if it is possible to give a simpler description of the field $L_{G}$ that works for every $G$.

\section*{Acknowledgements}
In the first place we would like to thank the referee for the various constructive and
quite valuable comments and suggestions on the paper. We also thank David Dynerman, Evan Dummit and Jordan Ellenberg for
helpful comments on an earlier version of this paper, and the organizers of the group theory seminar at UW-Madison for allowing us to present the results of this paper, and for their helpful
feedback.

\end{document}